\documentclass{amsart}

\newtheorem{theorem}{Theorem}[section]

\newtheorem{corollary}[theorem]{Corollary}

\theoremstyle{definition}

\theoremstyle{remark}

\numberwithin{equation}{section}

\begin{document}
\title{Identities of symmetry for $q$-Euler polynomials}{$\begin{array}{c}
         \text{}\\
         \text{}\
       \end{array}$}

\author{dae san kim} \thanks{}
\address{Department of Mathematics, Sogang University, Seoul 121-742, Korea}
\email{dskim@sogong.ac.kr}
\address{}
\email{}
\thanks{}

\subjclass[]{}

\date{}

\dedicatory{}


\begin{abstract}
In this paper, we derive eight basic identities of symmetry in three variables related to
$q$-Euler polynomials and the $q$-analogue of alternating power sums. These and most of their corollaries are new, since there have been results only about identities of symmetry in two variables. These abundance of symmetries shed new light even on the existing identities so as to yield some further interesting ones. The derivations of identities are based on the
$p$-adic integral expression of the generating function for the $q$-Euler polynomials and the quotient of integrals that can be expressed as the exponential generating function for the
$q$-analogue of alternating power sums.

\end{abstract}

\maketitle

\section{Introduction and preliminaries}
 Let $p$ be a fixed odd prime. Throughout this paper,
$\mathbb{Z}_{p}$, $\mathbb{Q}_{p}$, $\mathbb{C}_{p}$ will respectively denote the ring of
$p$-adic integers, the field of $p$-adic rational numbers and the completion of the algebraic closure of
$\mathbb{Q}_{p}$. For a continuous function
$f:\mathbb{Z}_{p}\longrightarrow\mathbb{C}_{p}$,
the $p$-adic fermionic integral of
 $f$ is defined by
\begin{align*}
\int_{\mathbb{Z}_{p}}f(z)d\mu_{-1}(z)=\lim_{N\rightarrow\infty}
\sum^{p^{N}-1}_{j=0}f(j)(-1)^{j}.
\end{align*}
Then it is easy to see that

\begin{equation}\label{a1}
\int_{\mathbb{Z}_{p}}f(z+1)d\mu_{-1}(z)+\int_{\mathbb{Z}_{p}}f(z)d\mu_{-1}(z)=
2f(0).
\end{equation}

Let $|~|_p$ be the normalized absolute value of $\mathbb{C}_{p}$,
such that $|p|_{p}=\frac{1}{p}$, and let

\begin{equation}\label{a2}
E=\{t\in\mathbb{C}_{p}||t|_p< p^{\frac{-1}{p-1}}\}.
\end{equation}
Assume that $q,t\in\mathbb{C}_{p}$, with $q-1,t\in E$, so that
$q^z=\exp(z\log q)$ and $e^{zt}$ are, as functions of $z$, analytic functions on $\mathbb{Z}_{p}$.
By applying (\ref{a1}) to
$f$, with $f(z)=q^{z} e^{tz}$, we get the $p$-adic
integral expression of the generating function for
$q$-Euler numbers $E_{n,q}$:

\begin{equation}\label{a3}
\int_{\mathbb{Z}_{p}}q^{z}e^{zt}d\mu_{-1}(z)
=\frac{2}{qe^{t}+1}
=\sum^{\infty}_{n=0}E_{n,q}\frac{t^{n}}{n!}~(t\in E).
\end{equation}

So we have the following $p$-adic integral
expression of the generating function for the $q$-Euler polynomials $E_{n,q}(x)$:

\begin{equation}\label{a4}
\int_{\mathbb{Z}_{p}}q^{z}e^{(x+z)t}d\mu_{-1}(z)
=\frac{2}{qe^{t}+1}e^{xt}
=\sum^{\infty}_{n=0}E_{n,q}(x)\frac{t^{n}}{n!}~(t\in E, x\in\mathbb{Z}_{p}).
\end{equation}
Note here that in \cite{A7} $\zeta$ was used in place of
$q$, and that $q$-Euler numbers and polynomials were coined respectively as
$\zeta$-Euler numbers and polynomials.

  Let $T_{k,q}(n)$ denote the $q$-analogue of alternating $k$th power sum of the first
$n+1$ nonnegative integers, namely

\begin{equation}\label{a5}
T_{k,q}(n)=\sum_{i=0}^{n}(-1)^{i}i^{k}q^{i}=(-1)^{0}0^{k}q^{0}+(-1)^{1}1^{k}q^{1}+\cdots +(-1)^{n}n^{k}q^{n}.
\end{equation}
In particular,
\begin{equation}\label{a6}
T_{0,q}(n)=\frac{(-q)^{n+1}-1}{(-q)-1}=[n+1]_{-q},~~
T_{k,q}(0)=
\begin{cases}
\begin{split}
1,& \quad \text {for}~ k=0,\\
0,& \quad \text {for}~ k>0.
\end{split}
\end{cases}
\end{equation}
From (\ref{a3}) and (\ref{a5}), one easily derives the following identities: for any odd positive integer
$w$,
\begin{equation}\label{a7}
\frac{\int_{\mathbb{Z}_{p}}q^{x}e^{xt}d\mu_{-1}(x)}{\int_{\mathbb{Z}_{p}}
q^{wy}e^{wyt}d\mu_{-1}(y)}=
\sum_{i=0}^{w-1}(-1)^{i}q^{i}e^{it}=\sum_{k=0}^{\infty}T_{k,q}(w-1)\frac{t^{k}}{k!}~
(t\in E).
\end{equation}
In what follows, we will always assume that the $p$-adic fermionic
 integrals of the various exponential functions on $\mathbb{Z}_{p}$ are defined for
$t\in E$(cf. (\ref{a2})), and therefore it will not be mentioned.

  \cite{A1}, \cite{A2}, \cite{A5}, \cite{A8} and \cite{A9} are some of the previous works on identities of symmetry
  involving Bernoulli polynomials and power sums. These results were generalized in
\cite{A4} to obtain identities of symmetry involving three variables in contrast to the previous works involving just two variables.

  In this paper, we will produce  8 basic identities of symmetry in three variables $w_1$,$w_2$,$w_3$
  related to $q$-Euler polynomials and the $q$-analogue of alternating power sums
  (cf. (\ref{a44}), (\ref{a45}), (\ref{a48}), (\ref{a51}), (\ref{a55}), (\ref{a57}),
  (\ref{a59}), (\ref{a60})). These and most of their corollaries seem to be new,  since there have been results only about identities of symmetry in two variables in the literature. These abundance of symmetries shed new light even on the existing identities. For instance, it has been known that (\ref{a8}) and (\ref{a9}) are equal and (\ref{a10}) and (\ref{a11}) are so(cf. [7, (2.11),(2,16)]).
In fact, (\ref{a8})-(\ref{a11}) are all equal, as they can be derived from one and the same $p$-adic integral.
Perhaps, this was neglected to mention in \cite{A7}. In addition, we have a bunch of new identities in (\ref{a12})-(\ref{a15}).
All of these were obtained as corollaries(cf. Cor. \ref{b9}, \ref{b12}, \ref{b15}) to some of the basic identities by specializing the variable $w_3$ as 1. Those would not be unearthed if more symmetries had not been available.  Related to
$q$-Bernoulli polynomials and the $q$-analogue of power sums,  identities of symmetry in three variables were also obtained in
\cite{A3} as an extension of identities of symmetry in two variables in
\cite{A6}.

Let $w_{1},w_{2}$ be any odd positive integers.  Then we have:

\begin{align}
\label{a8}
&\sum_{k=0}^{n}\binom{n}{k}E_{k,q^{w_{2}}}(w_{1}y_{1})T_{n-k,q^{w_{1}}}(w_{2}-1)w_{1}^{n-k}w_{2}^{k}\\
\label{a9}
=&\sum_{k=0}^{n}\binom{n}{k}E_{k,q^{w_{1}}}(w_{2}y_{1})T_{n-k,q^{w_{2}}}(w_{1}-1)w_{2}^{n-k}w_{1}^{k}\\
\label{a10}
=&w_{1}^{n}\sum_{i=0}^{w_{1}-1}(-1)^{i}q^{w_{2}i}E_{n,q^{w_{1}}}(w_{2}y_{1}+\frac{w_{2}}{w_{1}}i)\\
\label{a11}
=&w_{2}^{n}\sum_{i=0}^{w_{2}-1}(-1)^{i}q^{w_{1}i}E_{n,q^{w_{2}}}(w_{1}y_{1}+\frac{w_{1}}{w_{2}}i)\\
\label{a12}
=&\sum_{k+l+m=n}^{}\binom{n}{k,l,m}E_{k,q^{w_{1}w_{2}}}(y_{1})T_{l,q^{w_{2}}}(w_{1}-1)T_{m,q^{w_{1}}}(w_{2}-1)w_{1}^{k+m}w_{2}^{k+l}\\
\label{a13}
=&w_{1}^{n}\sum_{k=0}^{n}\binom{n}{k}T_{n-k,q^{w_{1}}}(w_{2}-1)w_{2}^{k}\sum_{i=0}^{w_{1}-1}(-1)^{i}q^{w_{2}i}E_{k,q^{w_{1}w_{2}}}
(y_{1}+\frac{i}{w_{1}})\\
\label{a14}
=&w_{2}^{n}\sum_{k=0}^{n}\binom{n}{k}T_{n-k,q^{w_{2}}}(w_{1}-1)w_{1}^{k}\sum_{i=0}^{w_{2}-1}(-1)^{i}q^{w_{1}i}E_{k,q^{w_{1}w_{2}}}
(y_{1}+\frac{i}{w_{2}})\\
\label{a15}
=&(w_{1}w_{2})^{n}\sum_{i=0}^{w_{1}-1}\sum_{j=0}^{w_{2}-1}(-1)^{i+j}q^{w_{2}i+w_{1}j}E_{n,q^{w_{1}w_{2}}}(y_{1}+\frac{i}{w_{1}}+\frac{j}{w_{2}}).
\end{align}

The derivations of identities are based on the $p$-adic integral
expression of the generating function for the $q$-Euler polynomials
in (\ref{a4}) and the quotient of integrals in (\ref{a7}) that can
be expressed as the exponential generating function for the
$q$-analogue of alternating power sums. We indebted this idea to the
papers [5,6].

\section{Several types of quotients of fermionic integrals}

  Here we will introduce several types of quotients of $p$-adic fermionic
 integrals on $\mathbb{Z}_{p}$ or $\mathbb{Z}_{p}^{3}$ from which some
interesting identities follow owing to the built-in symmetries in
$w_{1},w_{2},w_{3}$. In the following, $w_{1},w_{2},w_{3}$ are all
positive integers and all of the explicit expressions of integrals
in (\ref{a17}), (\ref{a19}), (\ref{a21}), and (\ref{a23}) are
obtained from the identity in (\ref{a3}).
\\
\\
(a) Type $\Lambda_{23}^{i}$ (for $i=0,1,2,3$)\\
$I(\Lambda_{23}^{i})$
\begin{align}\label{a16}
&=\frac{\begin{array}{c}
          \int_{\mathbb
{Z}_{p}^{3}}q^{w_{2}w_{3}x_{1}+w_{1}w_{3}x_{2}+w_{1}w_{2}x_{3}}e^{(w_{2}w_{3}x_{1}
+w_{1}w_{3}x_{2}+w_{1}w_{2}x_{3}
+w_{1}w_{2}w_{3}(\sum_{j=1}^{3-i}y_{j}))t} \\
          \qquad\qquad\qquad\qquad\qquad\qquad\qquad\qquad d\mu_{-1}(x_{1})d\mu_{-1}(x_{2})d\mu_{-1}(x_{3})
        \end{array}
}
 {(\int_{\mathbb
{Z}_{p}}q^{w_{1}w_{2}w_{3}x_{4}}e^{w_{1}w_{2}w_{3}x_{4}t}d\mu_{-1}(x_{4}))^{i}}
\end{align}

\begin{align}
\label{a17}
=\frac{2^{3-i}e^{w_{1}w_{2}w_{3}(\sum_{j=1}^{3-i}y_{j})t}(q^{w_{1}w_{2}w_{3}}
e^{w_{1}w_{2}w_{3}t}+1)^{i}}
{(q^{w_{2}w_{3}}e^{w_{2}w_{3}t}+1)(q^{w_{1}w_{3}}e^{w_{1}w_{3}t}+1)
(q^{w_{1}w_{2}}e^{w_{1}w_{2}t}+1)}.\qquad\qquad\qquad\quad
\end{align}
\\
\\
(b) Type $\Lambda_{13}^{i}$ (for $i=0,1,2,3$)\\
$I(\Lambda_{13}^{i})$
\begin{align}
\label{a18} &=\frac{\begin{array}{c}
                      \int_{\mathbb
{Z}_{p}^{3}}q^{w_{1}x_{1}+w_{2}x_{2}+w_{3}x_{3}}e^{(w_{1}x_{1}+w_{2}x_{2}
+w_{3}x_{3}+w_{1}w_{2}w_{3} (\sum_{j=1}^{3-i}y_{j}))t} \\
                      \qquad\qquad\qquad\qquad\qquad\qquad d\mu_{-1}
(x_{1})d\mu_{-1} (x_{2})d\mu_{-1} (x_{3})
                    \end{array}
}
{(\int_{\mathbb{Z}_{p}}q^{w_{1}w_{2}w_{3}x_{4}}
e^{w_{1}w_{2}w_{3}x_{4}t}d\mu_{-1} (x_{4}))^{i}}\\
\label{a19}
&=\frac{2^{3-i}e^{w_{1}w_{2}w_{3}(\sum_{j=1}^{3-i}y_{j})t}(q^{w_{1}w_{2}w_{3}}
e^{w_{1}w_{2}w_{3}t}+1)^{i}}
{(q^{w_{1}}e^{w_{1}t}+1)(q^{w_{2}}e^{w_{2}t}+1)(q^{w_{3}}e^{w_{3}t}+1)}.
\end{align}
\\
\\
(c-0) Type $\Lambda_{12}^{0}$\\
$I(\Lambda_{12}^{0})$
\begin{equation}
\begin{split}
\label{a20} =\int_{\mathbb
{Z}_{p}^{3}}&q^{w_{1}x_{1}+w_{2}x_{2}+w_{3}x_{3}}\\
\times&e^{(w_{1}x_{1}+w_{2}x_{2}+w_{3}x_{3}+w_{2}w_{3}y+w_{1}w_{3}y+w_{1}w_{2}y)t}
d\mu_{-1}(x_{1})d\mu_{-1}(x_{2})d\mu_{-1}(x_{3})
\end{split}
\end{equation}

\begin{equation}
\label{a21} =\frac{8e^{({w_{2}w_{3}+w_{1}w_{3}+w_{1}w_{2}})yt}}{
(q^{w_{1}}e^{w_{1}t}+1)(q^{w_{2}}e^{w_{2}t}+1)(q^{w_{3}}e^{w_{3}t}+1)}.\qquad\qquad\qquad\qquad\qquad\qquad\quad
\end{equation}
\\
\\
(c-1) Type $\Lambda_{12}^{1}$\\
$I(\Lambda_{12}^{1})$
\begin{align}
\label{a22}
&=\frac{\int_{\mathbb
{Z}_{p}^{3}}q^{w_{1}x_{1}+w_{2}x_{2}+w_{3}x_{3}}e^{
(w_{1}x_{1}+w_{2}x_{2}+w_{3}x_{3})t}d\mu_{-1}(x_{1})
d\mu_{-1}(x_{2})d\mu_{-1}(x_{3})} {\int_{\mathbb
{Z}_{p}^{3}}q^{w_{2}w_{3}z_{1}+w_{1}w_{3}z_{2}+w_{1}w_{2}z_{3}}e^{
(w_{2}w_{3}z_{1}+w_{1}w_{3}z_{2}+w_{1}w_{2}z_{3})t}
d\mu_{-1}(z_{1})d\mu_{-1}(z_{2})d\mu_{-1}(z_{3})}
\\
\label{a23} &=\frac{(q^{w_{2}w_{3}}e^{w_{2}w_{3}t}+1)
(q^{w_{1}w_{3}}e^{w_{1}w_{3}t}+1)(q^{w_{1}w_{2}}e^{w_{1}w_{2}t}+1)}
{(q^{w_{1}}e^{w_{1}t}+1)(q^{w_{2}}e^{w_{2}t}+1)(q^{w_{3}}e^{w_{3}t}+1)}.
\end{align}

  All of the above $p$-adic integrals of various types are invariant under
all permutations of $w_{1},w_{2},w_{3}$, as one can see either from
$p$-adic integral representations in (\ref{a16}), (\ref{a18}),
(\ref{a20}), and (\ref{a22}) or from their explicit evaluations in
(\ref{a17}), (\ref{a19}), (\ref{a21}), and (\ref{a23}).

\section{Identities for $q$-Euler polynomials}

 In the following $w_{1},w_{2},w_{3}$ are all odd positive integers except for
$(a-0)$ and $(c-0)$, where they are any positive integers.

$(a-0)$First, let's consider Type $\Lambda_{23}^{i}$, for each
$i=0,1,2,3.$ The following results can be easily obtained from
(\ref{a4}) and (\ref{a7}).
\begin{equation*}
\begin{split}
I(\Lambda_{23}^{0})=\int_{\mathbb
{Z}_{p}}q^{w_{2}w_{3}x_{1}}e^{w_{2}w_{3}(x_{1}+w_{1}y_{1})t}d\mu_{-1}(x_{1})
&\int_{\mathbb
{Z}_{p}}q^{w_{1}w_{3}x_{2}}e^{w_{1}w_{3}(x_{2}+w_{2}y_{2})t}d\mu_{-1}(x_{2})
\quad\qquad\qquad\qquad\\
\times&\int_{\mathbb Z_{p}}q^{w_{1}w_{2}x_{3}}e^{w_{1}w_{2}(x_{3}+w_{3}y_{3})t}
d\mu_{-1}(x_{3})\\
\end{split}
\end{equation*}
\begin{equation*}
\begin{split}
=(\sum_{k=0}^{\infty}\frac{E_{k,q^{w_{2}w_{3}}}(w_{1}y_{1})}{k!}
(w_{2}w_{3}t)^k)
(\sum_{l=0}^{\infty}&\frac{E_{l,q^{w_{1}w_{3}}}(w_{2}y_{2})}{l!}
(w_{1}w_{3}t)^l)\\
\times&(\sum_{m=0}^{\infty}\frac{E_{m,q^{w_{1}w_{2}}}(w_{3}y_{3})}{m!}
(w_{1}w_{2}t)^m)\\
\end{split}
\end{equation*}
\begin{equation}\label{a24}
\begin{split}
=\sum_{n=0}^{\infty}(\sum_{k+l+m=n}\binom{n}{k,l,m}
E_{k,q^{w_{2}w_{3}}}&(w_{1}y_{1})E_{l,q^{w_{1}w_{3}}}(w_{2}y_{2})\\
\times&E_{m,q^{w_{1}w_{2}}}(w_{3}y_{3})w_{1}^{l+m}w_{2}^{k+m}w_{3}^{k+l})
\frac{t^{n}}{n!},~\quad\\
\end{split}
\end{equation}
where the inner sum is over all nonnegative integers $k,l,m$, with
$k+l+m=n$, and
\begin{equation}\label{a25}
\binom{n}{k,l,m}=\frac{n!}{k!l!m!}.
\end{equation}
\\
\\
(a-1) Here we write $I(\Lambda_{23}^{1})$ in two different ways:
\\
\\
(1) $I(\Lambda_{23}^{1})$
\begin{equation}\label{a26}
\begin{split}
\qquad=\int_{\mathbb
{Z}_{p}}q^{w_{2}w_{3}x_{1}}e^{w_{2}w_{3}(x_{1}+w_{1}y_{1})t}d\mu_{-1}(x_{1})
&\int_{\mathbb{Z}_{p}}q^{w_{1}w_{3}x_{2}}e^{w_{1}w_{3}(x_{2}+w_{2}y_{2})t}
d\mu_{-1}(x_{2})\quad\qquad\qquad\qquad\\
\times&\frac{\int_{\mathbb
{Z}_{p}}q^{w_{1}w_{2}x_{3}}e^{w_{1}w_{2}x_{3}t}d\mu_{-1}(x_{3})}
{\int_{\mathbb{Z}_{p}}q^{w_{1}w_{2}w_{3}x_{4}}
e^{w_{1}w_{2}w_{3}x_{4}t}d\mu_{-1}(x_{4})}
\end{split}
\end{equation}
\begin{equation*}
\begin{split}
=(\sum_{k=0}^{\infty}{E_{k,q^{w_{2}w_{3}}}(w_{1}y_{1})\frac{(w_{2}w_{3}t)^k}{k!}})
&(\sum_{l=0}^{\infty}{E_{l,q^{w_{1}w_{3}}}(w_{2}y_{2})\frac{(w_{1}w_{3}t)^l}{l!}})\\
\times&(\sum_{m=0}^{\infty}{T_{m,q^{w_{1}w_{2}}}(w_{3}-1)\frac{(w_{1}w_{2}t)^m}{m!}})\quad\quad\quad
\end{split}
\end{equation*}
\begin{equation}\label{a27}
\begin{split}
\qquad=\sum_{n=0}^{\infty}(\sum_{k+l+m=n}\binom{n}
{k,l,m}&E_{k,q^{w_{2}w_{3}}}(w_{1}y_{1})E_{l,q^{w_{1}w_{3}}}(w_{2}y_{2})\\
&\times T_{m,q^{w_{1}w_{2}}}(w_{3}-1)w_{1}^{l+m}w_{2}^{k+m}w_{3}^{k+l})
\frac{t^{n}}{n!}.\qquad\qquad
\end{split}
\end{equation}
\\
\\
(2) Invoking (\ref{a7}), (\ref{a26}) can also be written as
\\
\\
$I(\Lambda_{23}^{1})$
\begin{equation}\label{a28}
\begin{split}
=\sum_{i=0}^{w_{3}-1}(-1)^{i}q^{w_{1}w_{2}i}
\int_{\mathbb
{Z}_{p}}q^{w_{2}w_{3}x_{1}}&e^{w_{2}w_{3}(x_{1}+w_{1}y_{1})t}d\mu_{-1}(x_{1})\qquad\qquad\\
&\times\int_{\mathbb{Z}_{p}}q^{w_{1}w_{3}x_{2}}e^{w_{1}w_{3}(x_{2}+w_{2}y_{2}+\frac{w_{2}}{w_{3}}i)t}d\mu_{-1}(x_{2})\qquad\qquad\qquad\\
\end{split}
\end{equation}
\begin{equation*}
\begin{split}
=\sum_{i=0}^{w_{3}-1}(-1)^{i}q^{w_{1}w_{2}i}
(\sum_{k=0}^{\infty}E_{k,q^{w_{2}w_{3}}}&(w_{1}y_{1})\frac{(w_{2}w_{3}t)^k}{k!})\\
\times&(\sum_{l=0}^{\infty}E_{l,q^{w_{1}w_{3}}}(w_{2}y_{2}+\frac{w_{2}}{w_{3}}i)\frac{(w_{1}w_{3}t)^l}{l!})\\
\end{split}
\end{equation*}
\begin{equation*}
\begin{split}
=\sum_{n=0}^{\infty}(w_{3}^{n}\sum_{k=0}^{n}\binom{n}{k}E_{k,q^{w_{2}w_{3}}}&(w_{1}y_{1})
\sum_{i=0}^{w_{3}-1}(-1)^{i}q^{w_{1}w_{2}i}\\
\times&E_{n-k,q^{w_{1}w_{3}}}(w_{2}y_{2}+\frac{w_{2}}{w_{3}}i)w_{1}^{n-k}w_{2}^{k})
\frac{t^{n}}{n!}.\qquad\\
\end{split}
\end{equation*}
(a-2) Here we write $I(\Lambda_{23}^{2})$ in three different ways:
\\
\\
(1) $I(\Lambda_{23}^{2})$
\begin{equation}\label{a29}
\begin{split}
=\int_{\mathbb
{Z}_{p}}&q^{w_{2}w_{3}x_{1}}e^{w_{2}w_{3}(x_{1}+w_{1}y_{1})t}d\mu_{-1}(x_{1})\\
\times&\frac{\int_{\mathbb
{Z}_{p}}q^{w_{1}w_{3}x_{2}}e^{w_{1}w_{3}x_{2}t}d\mu_{-1}(x_{2})}{\int_{\mathbb
{Z}_{p}}q^{w_{1}w_{2}w_{3}x_{4}}e^{w_{1}w_{2}w_{3}x_{4}t}d\mu_{-1}(x_{4})} \times
\frac{\int_{\mathbb
{Z}_{p}}q^{w_{1}w_{2}x_{3}}e^{w_{1}w_{2}x_{3}t}d\mu_{-1}(x_{3})}{\int_{\mathbb
{Z}_{p}}q^{w_{1}w_{2}w_{3}x_{4}}e^{w_{1}w_{2}w_{3}x_{4}t}d\mu_{-1}(x_{4})}
\end{split}
\end{equation}
\begin{equation*}
\begin{split}
\quad=(\sum_{k=0}^{\infty}&{E_{k,q^{w_{2}w_{3}}}(w_{1}y_{1})\frac{(w_{2}w_{3}t)^k}{k!}})\\
\times&(\sum_{l=0}^{\infty}{T_{l,q^{w_{1}w_{3}}}(w_{2}-1)\frac{(w_{1}w_{3}t)^l}{l!}})
(\sum_{m=0}^{\infty}{T_{m,q^{w_{1}w_{2}}}(w_{3}-1)\frac{(w_{1}w_{2}t)^m}{m!}})\qquad\quad
\end{split}
\end{equation*}
\begin{equation}\label{a30}
\begin{split}
=\sum_{n=0}^{\infty}(\sum_{k+l+m=n}\binom{n}{k,l,m}E_{k,q^{w_{2}w_{3}}}(w_{1}y_{1})T_{l,q^{w_{1}w_{3}}}&(w_{2}-1)
T_{m,q^{w_{1}w_{2}}}(w_{3}-1)\\
\times&w_{1}^{l+m}w_{2}^{k+m}w_{3}^{k+l})\frac{t^{n}}{n!}.\qquad\qquad
\end{split}
\end{equation}
\\
\\

(2) Invoking (\ref{a7}), (\ref{a29}) can also be written as
\\
\\
~$I(\Lambda_{23}^{2})$
\begin{equation}\label{a31}
\begin{split}
=\sum_{i=0}^{w_{2}-1}(-1)^{i}q^{w_{1}w_{3}i}\int_{\mathbb
{Z}_{p}}q^{w_{2}w_{3}x_{1}}&e^{w_{2}w_{3}(x_{1}+w_{1}y_{1}+\frac{w_{1}}{w_{2}}i)t}d\mu_{-1}(x_{1})\\
&\times\frac{\int_{\mathbb
{Z}_{p}}q^{w_{1}w_{2}x_{3}}e^{w_{1}w_{2}x_{3}t}d\mu_{-1}(x_{3})}{\int_{\mathbb
{Z}_{p}}q^{w_{1}w_{2}w_{3}x_{4}}e^{w_{1}w_{2}w_{3}x_{4}t}d\mu_{-1}(x_{4})}\quad
\end{split}
\end{equation}
\begin{equation*}
\begin{split}
=\sum_{i=0}^{w_{2}-1}(-1)^{i}q^{w_{1}w_{3}i}(\sum_{k=0}^{\infty}E_{k,q^{w_{2}w_{3}}}&(w_{1}y_{1}+\frac{w_{1}}{w_{2}}i)
\frac{(w_{2}w_{3}t)^k}{k!})\\
&\times(\sum_{l=0}^{\infty}T_{l,q^{w_{1}w_{2}}}(w_{3}-1)\frac{(w_{1}w_{2}t)^l}{l!})
\end{split}
\end{equation*}
\begin{equation}\label{a32}
\begin{split}
=\sum_{n=0}^{\infty}&(w_{2}^{n}\sum_{k=0}^{n}\binom{n}{k}\\
&\times\sum_{i=0}^{w_{2}-1}(-1)^{i}q^{w_{1}w_{3}i}E_{k,q^{w_{2}w_{3}}}(w_{1}y_{1}+\frac{w_{1}}{w_{2}}i)
T_{n-k,q^{w_{1}w_{2}}}(w_{3}-1)w_{1}^{n-k}w_{3}^{k})\frac{t^{n}}{n!}.
\end{split}
\end{equation}
\\
\\
(3) Invoking (\ref{a7}) once again, (\ref{a31}) can be written as
\\
\\
~$I(\Lambda_{23}^{2})$
\begin{equation*}
=\sum_{i=0}^{w_{2}-1}\sum_{j=0}^{w_{3}-1}(-1)^{i+j}q^{w_{1}(w_{3}i+w_{2}j)}
\int_{\mathbb
{Z}_{p}}q^{w_{2}w_{3}x_{1}}e^{w_{2}w_{3}(x_{1}+w_{1}y_{1}+\frac{w_{1}}{w_{2}}i+\frac{w_{1}}{w_{3}}j)t}d\mu_{-1}(x_{1})
\qquad\qquad
\end{equation*}
\begin{equation*}
=\sum_{i=0}^{w_{2}-1}\sum_{j=0}^{w_{3}-1}(-1)^{i+j}q^{w_{1}(w_{3}i+w_{2}j)}
\sum_{n=0}^{\infty}E_{n,q^{w_{2}w_{3}}}(w_{1}y_{1}+\frac{w_{1}}{w_{2}}i+\frac{w_{1}}{w_{3}}j)\frac{(w_{2}w_{3}t)^n}{n!})
\qquad\qquad\qquad\qquad
\end{equation*}
\begin{equation}\label{a33}
=\sum_{n=0}^{\infty}((w_{2}w_{3})^{n}\sum_{i=0}^{w_{2}-1}\sum_{j=0}^{w_{3}-1}(-1)^{i+j}
q^{w_{1}(w_{3}i+w_{2}j)}
E_{n,q^{w_{2}w_{3}}}(w_{1}y_{1}+
\frac{w_{1}}{w_{2}}i+\frac{w_{1}}{w_{3}}j))\frac{t^{n}}{n!}.\qquad
\end{equation}
\\
\\
(a-3)
\\
\\
~$I(\Lambda_{23}^{3})$
\begin{equation*}
\begin{split}
=&\frac{\int_{\mathbb
{Z}_{p}}q^{w_{2}w_{3}x_{1}}e^{w_{2}w_{3}x_{1}t}d\mu_{-1}(x_{1})}
{\int_{\mathbb
{Z}_{p}}q^{w_{1}w_{2}w_{3}x_{4}}e^{w_{1}w_{2}w_{3}x_{4}t}d\mu_{-1}(x_{4})}\\
&\times\frac{\int_{\mathbb
{Z}_{p}}q^{w_{1}w_{3}x_{2}}e^{w_{1}w_{3}x_{2}t}d\mu_{-1}(x_{2})}{\int_{\mathbb
{Z}_{p}}q^{w_{1}w_{2}w_{3}x_{4}}e^{w_{1}w_{2}w_{3}x_{4}t}d\mu_{-1}(x_{4})}
\times\frac{\int_{\mathbb
{Z}_{p}}q^{w_{1}w_{2}x_{3}}e^{w_{1}w_{2}x_{3}t}d\mu_{-1}(x_{3})}{\int_{\mathbb
{Z}_{p}}q^{w_{1}w_{2}w_{3}x_{4}}e^{w_{1}w_{2}w_{3}x_{4}t}d\mu_{-1}(x_{4})}
\end{split}
\end{equation*}
\begin{equation*}
\begin{split}
=&(\sum_{k=0}^{\infty}T_{k,q^{w_{2}w_{3}}}(w_{1}-1)\frac{(w_{2}w_{3}t)^{k}}{k!})\\
\times&(\sum_{l=0}^{\infty}T_{l,q^{w_{1}w_{3}}}(w_{2}-1)\frac{(w_{1}w_{3}t)^{l}}{l!})
(\sum_{m=0}^{\infty}T_{m,q^{w_{1}w_{2}}}(w_{3}-1)\frac{(w_{1}w_{2}t)^{m}}{m!})\qquad\qquad
\end{split}
\end{equation*}
\begin{equation}\label{a34}
\begin{split}
=\sum_{n=0}^{\infty}(\sum_{k+l+m=n}^{}\binom{n}{k,l,m}T_{k,q^{w_{2}w_{3}}}(w_{1}-1)T_{l,q^{w_{1}w_{3}}}&(w_{2}-1)
T_{m,q^{w_{1}w_{2}}}(w_{3}-1)\\
\times&w_{1}^{l+m}w_{2}^{k+m}w_{3}^{k+l})\frac{t^{n}}{n!}.\quad
\end{split}
\end{equation}
\\
\\
(b) For Type $\Lambda_{13}^{i}~(i=0,1,2,3)$, we may consider the
analogous things to the ones in (a-0), (a-1), (a-2), and (a-3).
However, these do not lead us to new identities. Indeed, if we
substitute $w_{2}w_{3},w_{1}w_{3},w_{1}w_{2}$ respectively for
$w_{1},w_{2},w_{3}$ in (\ref{a16}), this amounts to replacing $t$ by
$w_{1}w_{2}w_{3}t$ and $q$ by $q^{w_{1}w_{2}w_{3}}$ in (\ref{a18}). So, upon replacing
$w_{1},w_{2},w_{3}$ respectively by
$w_{2}w_{3},w_{1}w_{3},w_{1}w_{2}$, and then dividing by
$(w_{1}w_{2}w_{3})^n$ and replacing
$q^{w_{1}w_{2}w_{3}}$ by $q$, in each of the expressions of Theorem
\ref{b1} through Corollary \ref{b15},
we will get the corresponding symmetric identities for Type
$\Lambda_{13}^{i}~(i=0,1,2,3)$.
\\
\\
(c-0)
\\
$I(\Lambda_{12}^{0})$
\begin{align*}
\begin{split}
=\int_{\mathbb{Z}_{p}}q^{w_{1}x_{1}}&e^{w_{1}(x_{1}+w_{2}y)t}d\mu_{-1}(x_{1})\\
\times&\int_{\mathbb{Z}_{p}}q^{w_{2}x_{2}}e^{w_{2}(x_{2}+w_{3}y)t}d\mu_{-1}(x_{2})
\int_{\mathbb Z_{p}}q^{w_{3}x_{3}}e^{w_{3}(x_{3}+w_{1}y)t}d\mu_{-1}(x_{3})\qquad\qquad
\end{split}
\end{align*}
\begin{align*}
=(\sum_{k=0}^{\infty}\frac{E_{k,q^{w_{1}}}(w_{2}y)}{k!}(w_{1}t)^{k})
(\sum_{l=0}^{\infty}\frac{E_{l,q^{w_{2}}}(w_{3}y)}{l!}(w_{2}t)^{l})
(\sum_{m=0}^{\infty}\frac{E_{m,q^{w_{3}}}(w_{1}y)}{m!}(w_{3}t)^{m})\qquad\qquad
\end{align*}
\begin{align}\label{a35}
=\sum_{n=0}^{\infty}(\sum_{k+l+m=n}\binom{n}{k,l,m}E_{k,q^{w_{1}}}(w_{2}y)E_{l,q^{w_{2}}}
(w_{3}y)E_{m,q^{w_{3}}}(w_{1}y)w_{1}^{k}w_{2}^{l}w_{3}^{m})\frac{t^n}{n!}.~\qquad\qquad
\end{align}
\\
\\
(c-1)\\
\\
\begin{equation*}
\begin{split}
I(\Lambda_{12}^{1})=&\frac{\int_{\mathbb
Z_{p}}q^{w_{1}x_{1}}e^{w_{1}x_{1}t}d\mu_{-1}(x_{1})}{\int_{\mathbb
Z_{p}}q^{w_{1}w_{2}z_{3}}e^{w_{1}w_{2}z_{3}t}d\mu_{-1}(z_{3})}\\
\times&\frac{\int_{\mathbb
Z_{p}}q^{w_{2}x_{2}}e^{w_{2}x_{2}t}d\mu_{-1}(x_{2})}{\int_{\mathbb
Z_{p}}q^{w_{2}w_{3}z_{1}}e^{w_{2}w_{3}z_{1}t}d\mu_{-1}(z_{1})}
\times\frac{\int_{\mathbb
Z_{p}}q^{w_{3}x_{3}}e^{w_{3}x_{3}t}d\mu_{-1}(x_{3})}{\int_{\mathbb
Z_{p}}q^{w_{3}w_{1}z_{2}}e^{w_{3}w_{1}z_{2}t}d\mu_{-1}(z_{2})}\qquad\qquad
\end{split}
\end{equation*}
\begin{equation*}
\begin{split}
~\qquad=&(\sum_{k=0}^{\infty}T_{k,q^{w_{1}}}(w_{2}-1)\frac{(w_{1}t)^{k}}{k!})\\
\times&(\sum_{l=0}^{\infty}T_{l,q^{w_{2}}}(w_{3}-1)\frac{(w_{2}t)^{l}}{l!})
(\sum_{m=0}^{\infty}T_{m,q^{w_{3}}}(w_{1}-1)\frac{(w_{3}t)^{m}}{m!})\qquad\qquad\qquad
\end{split}
\end{equation*}
\begin{equation}\label{a36}
\begin{split}
=\sum_{n=0}^{\infty}(\sum_{k+l+m=n}^{}\binom{n}{k,l,m}&T_{k,q^{w_{1}}}(w_{2}-1)T_{l,q^{w_{2}}}(w_{3}-1)\\
&\times T_{m,q^{w_{3}}}(w_{1}-1)w_{1}^{k}w_{2}^{l}w_{3}^{m})\frac{t^{n}}{n!}.\qquad\qquad\qquad
\end{split}
\end{equation}
\section{Main theorems}
  As we noted earlier in the last paragraph of Section 2, the various
types of quotients of $p$-adic fermionic
 integrals are invariant
under any permutation of $w_{1},w_{2},w_{3}$. So the corresponding
expressions in Section 3 are also invariant under any permutation of
$w_{1},w_{2},w_{3}$. Thus our results about identities of symmetry
will be immediate consequences of this observation.

  However, not all permutations of an expression in Section 3 yield
distinct ones. In fact, as these expressions are obtained by
permuting $w_{1},w_{2},w_{3}$ in a single one labeled by them, they
can be viewed as a group in a natural manner and hence it is
isomorphic to a quotient of $S_{3}$. In particular, the number of
possible distinct expressions are $1,2,3,$ or $6$. (a-0), (a-1(1)),
(a-1(2)), and (a-2(2)) give the full six identities of symmetry,
(a-2(1)) and (a-2(3)) yield three identities of symmetry, and (c-0)
and (c-1) give two identities of symmetry, while the expression in
(a-3) yields no identities of symmetry.

  Here we will just consider the cases of Theorems \ref{b8} and \ref{b17}, leaving
the others as easy exercises for the reader. As for the case of
Theorem \ref{b8}, in addition to (\ref{a50})-(\ref{a52}), we get the
following three ones:
\begin{equation}\label{a37}
\begin{split}
\sum_{k+l+m=n}\binom{n}{k,l,m}E_{k,q^{w_{2}w_{3}}}(w_{1}y_{1})T_{l,q^{w_{1}w_{2}}}(w_{3}-1)&T_{m,q^{w_{1}w_{3}}}
(w_{2}-1)\\
&\times w_{1}^{l+m}w_{3}^{k+m}w_{2}^{k+l},
\end{split}
\end{equation}
\begin{equation}\label{a38}
\begin{split}
\sum_{k+l+m=n}\binom{n}{k,l,m}E_{k,q^{w_{1}w_{3}}}(w_{2}y_{1})T_{l,q^{w_{2}w_{3}}}(w_{1}-1)&T_{m,q^{w_{1}w_{2}}}
(w_{3}-1)\\
&\times w_{2}^{l+m}w_{1}^{k+m}w_{3}^{k+l},
\end{split}
\end{equation}
\begin{equation}\label{a39}
\begin{split}
\sum_{k+l+m=n}\binom{n}{k,l,m}E_{k,q^{w_{1}w_{2}}}(w_{3}y_{1})T_{l,q^{w_{1}w_{3}}}(w_{2}-1)
&T_{m,q^{w_{2}w_{3}}}(w_{1}-1)\\
&\times w_{3}^{l+m}w_{2}^{k+m}w_{1}^{k+l}.
\end{split}
\end{equation}
But, by interchanging $l$ and $m$, we see that (\ref{a37}),
(\ref{a38}), and (\ref{a39}) are respectively equal to (\ref{a50}),
(\ref{a51}), and (\ref{a52}).
\\
As to Theorem \ref{b17}, in addition to (\ref{a60}) and (\ref{a61}), we have:
\begin{align}
\label{a40}
&\sum_{k+l+m=n}\binom{n}{k,l,m}T_{k}(w_{2}-1)T_{l}(w_{3}-1)T_{m}
(w_{1}-1)w_{1}^{k}w_{2}^{l}w_{3}^{m},\\
\label{a41}
&\sum_{k+l+m=n}\binom{n}{k,l,m}T_{k}(w_{3}-1)T_{l}(w_{1}-1)T_{m}
(w_{2}-1)w_{2}^{k}w_{3}^{l}w_{1}^{m},\\
\label{a42}
&\sum_{k+l+m=n}\binom{n}{k,l,m}T_{k}(w_{3}-1)T_{l}(w_{2}-1)T_{m}
(w_{1}-1)w_{1}^{k}w_{3}^{l}w_{2}^{m},\\
\label{a43}
&\sum_{k+l+m=n}\binom{n}{k,l,m}T_{k}(w_{2}-1)T_{l}(w_{1}-1)T_{m}
(w_{3}-1)w_{3}^{k}w_{2}^{l}w_{1}^{m}.
\end{align}
\\
  However, (\ref{a40}) and (\ref{a41}) are equal to (\ref{a60}), as we
can see by applying the permutations $k\rightarrow l,l\rightarrow
m,m\rightarrow k$ for (\ref{a40}) and $k\rightarrow m,l\rightarrow
k,m\rightarrow l$ for (\ref{a41}). Similarly, we see that (\ref{a42})
and (\ref{a43}) are equal to (\ref{a61}), by applying permutations
$k\rightarrow l,l\rightarrow m,m\rightarrow k$ for (\ref{a42}) and
$k\rightarrow m,l\rightarrow k,m\rightarrow l$ for (\ref{a43}).
\begin{theorem}\label{b1}
Let $w_{1},w_{2},w_{3}$ be any positive integers. Then we have :

\begin{equation}\label{a44}
\begin{split}
&\sum_{k+l+m=n}\binom{n}{k,l,m}E_{k,q^{w_{2}w_{3}}}(w_{1}y_{1})E_{l,q^{w_{1}w_{3}}}
(w_{2}y_{2})E_{m,q^{w_{1}w_{2}}}(w_{3}y_{3})w_{1}^{l+m}w_{2}^{k+m}w_{3}^{k+l}\\
=&\sum_{k+l+m=n}\binom{n}{k,l,m}E_{k,q^{w_{2}w_{3}}}(w_{1}y_{1})E_{l,q^{w_{1}w_{2}}}
(w_{3}y_{2})E_{m,q^{w_{1}w_{3}}}(w_{2}y_{3})w_{1}^{l+m}w_{3}^{k+m}w_{2}^{k+l}\\
=&\sum_{k+l+m=n}\binom{n}{k,l,m}E_{k,q^{w_{1}w_{3}}}(w_{2}y_{1})E_{l,q^{w_{2}w_{3}}}
(w_{1}y_{2})E_{m,q^{w_{1}w_{2}}}(w_{3}y_{3})w_{2}^{l+m}w_{1}^{k+m}w_{3}^{k+l}\\
=&\sum_{k+l+m=n}\binom{n}{k,l,m}E_{k,q^{w_{1}w_{3}}}(w_{2}y_{1})E_{l,q^{w_{1}w_{2}}}
(w_{3}y_{2})E_{m,q^{w_{2}w_{3}}}(w_{1}y_{3})w_{2}^{l+m}w_{3}^{k+m}w_{1}^{k+l}\\
=&\sum_{k+l+m=n}\binom{n}{k,l,m}E_{k,q^{w_{1}w_{2}}}(w_{3}y_{1})E_{l,q^{w_{2}w_{3}}}
(w_{1}y_{2})E_{m,q^{w_{1}w_{3}}}(w_{2}y_{3})w_{3}^{l+m}w_{1}^{k+m}w_{2}^{k+l}\\
=&\sum_{k+l+m=n}\binom{n}{k,l,m}E_{k,q^{w_{1}w_{2}}}(w_{3}y_{1})E_{l,q^{w_{1}w_{3}}}
(w_{2}y_{2})E_{m,q^{w_{2}w_{3}}}(w_{1}y_{3})w_{3}^{l+m}w_{2}^{k+m}w_{1}^{k+l}.
\end{split}
\end{equation}
\end{theorem}
\begin{theorem}\label{b2}
  Let $w_{1},w_{2},w_{3}$ be any odd positive integers. Then we have:
\begin{equation}\label{a45}
\begin{split}
\sum_{k+l+m=n}\binom{n}{k,l,m}E_{k,q^{w_{2}w_{3}}}(w_{1}y_{1})E_{l,q^{w_{1}w_{3}}}
(w_{2}y_{2})&T_{m,q^{w_{1}w_{2}}}(w_{3}-1)\\
&\times w_{1}^{l+m}w_{2}^{k+m}w_{3}^{k+l}
\end{split}
\end{equation}
\begin{equation*}
\begin{split}
=\sum_{k+l+m=n}\binom{n}{k,l,m}E_{k,q^{w_{2}w_{3}}}(w_{1}y_{1})E_{l,q^{w_{1}w_{2}}}
(w_{3}y_{2})&T_{m,q^{w_{1}w_{3}}}(w_{2}-1)\\
&\times w_{1}^{l+m}w_{3}^{k+m}w_{2}^{k+l}
\end{split}
\end{equation*}
\begin{equation*}
\begin{split}
=\sum_{k+l+m=n}\binom{n}{k,l,m}E_{k,q^{w_{1}w_{3}}}(w_{2}y_{1})E_{l,q^{w_{2}w_{3}}}
(w_{1}y_{2})&T_{m,q^{w_{1}w_{2}}}(w_{3}-1)\\
&\times w_{2}^{l+m}w_{1}^{k+m}w_{3}^{k+l}
\end{split}
\end{equation*}
\begin{equation*}
\begin{split}
=\sum_{k+l+m=n}\binom{n}{k,l,m}E_{k,q^{w_{1}w_{3}}}(w_{2}y_{1})E_{l,q^{w_{1}w_{2}}}
(w_{3}y_{2})&T_{m,q^{w_{2}w_{3}}}(w_{1}-1)\\
&\times w_{2}^{l+m}w_{3}^{k+m}w_{1}^{k+l}
\end{split}
\end{equation*}
\begin{equation*}
\begin{split}
=\sum_{k+l+m=n}\binom{n}{k,l,m}E_{k,q^{w_{1}w_{2}}}(w_{3}y_{1})E_{l,q^{w_{1}w_{3}}}
(w_{2}y_{2})&T_{m,q^{w_{2}w_{3}}}(w_{1}-1)\\
&\times w_{3}^{l+m}w_{2}^{k+m}w_{1}^{k+l}\\
\end{split}
\end{equation*}
\begin{equation*}
\begin{split}
=\sum_{k+l+m=n}\binom{n}{k,l,m}E_{k,q^{w_{1}w_{2}}}(w_{3}y_{1})E_{l,q^{w_{2}w_{3}}}
(w_{1}y_{2})&T_{m,q^{w_{1}w_{3}}}(w_{2}-1)\\
&\times w_{3}^{l+m}w_{1}^{k+m}w_{2}^{k+l}.
\end{split}
\end{equation*}
\end{theorem}
Putting $w_{3}=1$ in (\ref{a45}), we get the following corollary.
\begin{corollary}\label{b3}
Let $w_{1},w_{2}$ be any odd positive integers.  Then we have:
\begin{equation}\label{a46}
\begin{split}
&\sum_{k=0}^{n}\binom{n}{k}E_{k,q^{w_{2}}}(w_{1}y_{1})E_{n-k,q^{w_{1}}}(w_{2}y_{2})w_{1}^{n-k}w_{2}^{k}\\
=&\sum_{k=0}^{n}\binom{n}{k}E_{k,q^{w_{1}}}(w_{2}y_{1})E_{n-k,q^{w_{2}}}(w_{1}y_{2})w_{2}^{n-k}w_{1}^{k}\\
=&\sum_{k+l+m=n}\binom{n}{k,l,m}E_{k,q^{w_{1}w_{2}}}(y_{1})E_{l,q^{w_{1}}}(w_{2}y_{2})
T_{m,q^{w_{2}}}(w_{1}-1)w_{2}^{k+m}w_{1}^{k+l}\\
=&\sum_{k+l+m=n}\binom{n}{k,l,m}E_{k,q^{w_{1}}}(w_{2}y_{1})E_{l,q^{w_{1}w_{2}}}(y_{2})
T_{m,q^{w_{2}}}(w_{1}-1)w_{2}^{l+m}w_{1}^{k+l}\\
=&\sum_{k+l+m=n}\binom{n}{k,l,m}E_{k,q^{w_{1}w_{2}}}(y_{1})E_{l,q^{w_{2}}}(w_{1}y_{2})
T_{m,q^{w_{1}}}(w_{2}-1)w_{1}^{k+m}w_{2}^{k+l}\\
=&\sum_{k+l+m=n}\binom{n}{k,l,m}E_{k,q^{w_{2}}}(w_{1}y_{1})E_{l,q^{w_{1}w_{2}}}(y_{2})
T_{m,q^{w_{1}}}(w_{2}-1)w_{1}^{l+m}w_{2}^{k+l}.
\end{split}
\end{equation}
\end{corollary}
Letting further $w_{2}=1$ in (\ref{a46}), we have the following
corollary.
\begin{corollary}\label{b4}
  Let $w_{1}$ be any odd positive integer.  Then we have:
\begin{equation}\label{a47}
\begin{split}
&\sum_{k=0}^{n}\binom{n}{k}E_{k,q}(w_{1}y_{1})E_{n-k,q^{w_{1}}}(y_{2})w_{1}^{n-k}\\
=&\sum_{k=0}^{n}\binom{n}{k}E_{k,q^{w_{1}}}(y_{1})E_{n-k,q}(w_{1}y_{2})w_{1}^{k}\\
=&\sum_{k+l+m=n}\binom{n}{k,l,m}E_{k,q^{w_{1}}}(y_{1})E_{l,q^{w_{1}}}(y_{2})T_{m,q}(w_{1}-1)w_{1}^{k+l}.
\end{split}
\end{equation}
\end{corollary}
\begin{theorem}\label{b5}
Let $w_{1},w_{2},w_{3}$ be any odd positive integers. Then we have:
\begin{equation}\label{a48}
\begin{split}
&w_{1}^{n}\sum_{k=0}^{n}\binom{n}{k}E_{k,q^{w_{1}w_{2}}}(w_{3}y_{1})w_{3}^{n-k}w_{2}^{k}
\sum_{i=0}^{w_{1}-1}(-1)^{i}q^{w_{2}w_{3}i}E_{n-k,q^{w_{1}w_{3}}}(w_{2}y_{2}+\frac{w_{2}}{w_{1}}i)\\
=&w_{1}^{n}\sum_{k=0}^{n}\binom{n}{k}E_{k,q^{w_{1}w_{3}}}(w_{2}y_{1})w_{2}^{n-k}w_{3}^{k}
\sum_{i=0}^{w_{1}-1}(-1)^{i}q^{w_{2}w_{3}i}E_{n-k,q^{w_{1}w_{2}}}(w_{3}y_{2}+\frac{w_{3}}{w_{1}}i)
\end{split}
\end{equation}
\begin{equation*}
\begin{split}
=&w_{2}^{n}\sum_{k=0}^{n}\binom{n}{k}E_{k,q^{w_{1}w_{2}}}(w_{3}y_{1})w_{3}^{n-k}w_{1}^{k}
\sum_{i=0}^{w_{2}-1}(-1)^{i}q^{w_{1}w_{3}i}E_{n-k,q^{w_{2}w_{3}}}(w_{1}y_{2}+\frac{w_{1}}{w_{2}}i)\\
=&w_{2}^{n}\sum_{k=0}^{n}\binom{n}{k}E_{k,q^{w_{2}w_{3}}}(w_{1}y_{1})w_{1}^{n-k}w_{3}^{k}
\sum_{i=0}^{w_{2}-1}(-1)^{i}q^{w_{1}w_{3}i}E_{n-k,q^{w_{1}w_{2}}}(w_{3}y_{2}+\frac{w_{3}}{w_{2}}i)\\
=&w_{3}^{n}\sum_{k=0}^{n}\binom{n}{k}E_{k,q^{w_{1}w_{3}}}(w_{2}y_{1})w_{2}^{n-k}w_{1}^{k}
\sum_{i=0}^{w_{3}-1}(-1)^{i}q^{w_{1}w_{2}i}E_{n-k,q^{w_{2}w_{3}}}(w_{1}y_{2}+\frac{w_{1}}{w_{3}}i)\\
=&w_{3}^{n}\sum_{k=0}^{n}\binom{n}{k}E_{k,q^{w_{2}w_{3}}}(w_{1}y_{1})w_{1}^{n-k}w_{2}^{k}
\sum_{i=0}^{w_{3}-1}(-1)^{i}q^{w_{1}w_{2}i}E_{n-k,q^{w_{1}w_{3}}}(w_{2}y_{2}+\frac{w_{2}}{w_{3}}i).
\end{split}
\end{equation*}
\end{theorem}
Letting $w_{3}=1$ in (\ref{a48}), we obtain alternative expressions
for the identities in (\ref{a46}).
\begin{corollary}\label{b6}
  Let $w_{1},w_{2}$ be any odd positive integers.  Then we have:
\begin{equation}\label{a49}
\begin{split}
&\sum_{k=0}^{n}\binom{n}{k}E_{k,q^{w_{2}}}(w_{1}y_{1})E_{n-k,q^{w_{1}}}(w_{2}y_{2})w_{1}^{n-k}w_{2}^{k}\\
=&\sum_{k=0}^{n}\binom{n}{k}E_{k,q^{w_{1}}}(w_{2}y_{1})E_{n-k,q^{w_{2}}}(w_{1}y_{2})w_{2}^{n-k}w_{1}^{k}\\
=&w_{1}^{n}\sum_{k=0}^{n}\binom{n}{k}E_{k,q^{w_{1}w_{2}}}(y_{1})w_{2}^{k}
\sum_{i=0}^{w_{1}-1}(-1)^{i}q^{w_{2}i}E_{n-k,q^{w_{1}}}(w_{2}y_{2}+\frac{w_{2}}{w_{1}}i)\\
=&w_{1}^{n}\sum_{k=0}^{n}\binom{n}{k}E_{k,q^{w_{1}}}(w_{2}y_{1})w_{2}^{n-k}
\sum_{i=0}^{w_{1}-1}(-1)^{i}q^{w_{2}i}E_{n-k,q^{w_{1}w_{2}}}(y_{2}+\frac{i}{w_{1}})\\
=&w_{2}^{n}\sum_{k=0}^{n}\binom{n}{k}E_{k,q^{w_{1}w_{2}}}(y_{1})w_{1}^{k}
\sum_{i=0}^{w_{2}-1}(-1)^{i}q^{w_{1}i}E_{n-k,q^{w_{2}}}(w_{1}y_{2}+\frac{w_{1}}{w_{2}}i)\\
=&w_{2}^{n}\sum_{k=0}^{n}\binom{n}{k}E_{k,q^{w_{2}}}(w_{1}y_{1})w_{1}^{n-k}
\sum_{i=0}^{w_{2}-1}(-1)^{i}q^{w_{1}i}E_{n-k,q^{w_{1}w_{2}}}(y_{2}+\frac{i}{w_{2}}).
\end{split}
\end{equation}
\end{corollary}
Putting further $w_{2}=1$ in (\ref{a49}), we have the alternative
expressions for the identities for (\ref{a47}).
\begin{corollary}\label{b7}
  Let $w_{1}$ be any odd positive integer.  Then we have:
\begin{equation*}
\begin{split}
&\sum_{k=0}^{n}\binom{n}{k}E_{k,q^{w_{1}}}(y_{1})E_{n-k,q}(w_{1}y_{2})w_{1}^{k}\\
=&\sum_{k=0}^{n}\binom{n}{k}E_{k,q^{w_{1}}}(y_{2})E_{n-k,q}(w_{1}y_{1})w_{1}^{k}\\
=&w_{1}^{n}\sum_{k=0}^{n}\binom{n}{k}E_{k,q^{w_{1}}}(y_{1})
\sum_{i=0}^{w_{1}-1}(-1)^{i}q^{i}E_{n-k,q^{w_{1}}}(y_{2}+\frac{i}{w_{1}}).
\end{split}
\end{equation*}
\end{corollary}
\begin{theorem}\label{b8}
Let $w_{1},w_{2},w_{3}$ be any odd positive integers. Then we have:
\begin{equation}
\begin{split}
\label{a50}
\sum_{k+l+m=n}\binom{n}{k,l,m}E_{k,q^{w_{2}w_{3}}}(w_{1}y_{1})&T_{l,q^{w_{1}w_{3}}}(w_{2}-1)
T_{m,q^{w_{1}w_{2}}}(w_{3}-1)\\
&\times w_{1}^{l+m}w_{2}^{k+m}w_{3}^{k+l}
\end{split}
\end{equation}
\begin{equation}
\begin{split}
\label{a51}
=\sum_{k+l+m=n}\binom{n}{k,l,m}E_{k,q^{w_{1}w_{3}}}(w_2y_{1})&T_{l,q^{w_{1}w_{2}}}(w_{3}-1)
T_{m,q^{w_{2}w_{3}}}(w_{1}-1)\\
&\times w_{2}^{l+m}w_{3}^{k+m}w_{1}^{k+l}
\end{split}
\end{equation}
\begin{equation}
\begin{split}
\label{a52}
=\sum_{k+l+m=n}\binom{n}{k,l,m}E_{k,q^{w_{1}w_{2}}}(w_{3}y_{1})&T_{l,q^{w_{2}w_{3}}}(w_{1}-1)
T_{m,q^{w_{1}w_{3}}}(w_{2}-1)\\
&\times w_{3}^{l+m}w_{1}^{k+m}w_{2}^{k+l}.
\end{split}
\end{equation}
\end{theorem}
Putting $w_{3}=1$ in (\ref{a50})-(\ref{a52}), we get the following
corollary.
\begin{corollary}\label{b9}
  Let $w_{1},w_{2}$ be any odd positive integers.  Then we have:
\begin{equation}\label{a53}
\begin{split}
&\sum_{k=0}^{n}\binom{n}{k}E_{k,q^{w_{2}}}(w_{1}y_{1})T_{n-k,q^{w_{1}}}(w_{2}-1)w_{1}^{n-k}w_{2}^{k}\\
=&\sum_{k=0}^{n}\binom{n}{k}E_{k,q^{w_{1}}}(w_{2}y_{1})T_{n-k,q^{w_{2}}}(w_{1}-1)w_{2}^{n-k}w_{1}^{k}\\
=&\sum_{k+l+m=n}\binom{n}{k,l,m}E_{k,q^{w_{1}w_{2}}}(y_{1})T_{l,q^{w_{2}}}(w_{1}-1)
T_{m,q^{w_{1}}}(w_{2}-1)w_{1}^{k+m}w_{2}^{k+l}.
\end{split}
\end{equation}
\end{corollary}
Letting further $w_{2}=1$ in (\ref{a53}), we get the following
corollary. This was also obtained in [7, (2.12)].
\begin{corollary}\label{b10}
  Let $w_{1}$ be any odd positive integer.  Then we have:
\begin{equation}\label{a54}
E_{n,q}(w_{1}y_{1})=\sum_{k=0}^{n}\binom{n}{k}E_{k,q^{w_{1}}}(y_{1})T_{n-k,q}(w_{1}-1)w_{1}^{k}.
\end{equation}
\end{corollary}
\begin{theorem}\label{b11}
Let $w_{1},w_{2},w_{3}$ be any odd positive integers. Then we have:
\begin{equation}\label{a55}
\begin{split}
w_{1}^{n}\sum_{k=0}^{n}\binom{n}{k}T_{n-k,q^{w_{1}w_{2}}}(w_{3}-1)w_{2}^{n-k}w_{3}^{k}
\sum_{i=0}^{w_{1}-1}(-1)^{i}q^{w_{2}w_{3}i}E_{k,q^{w_{1}w_{3}}}(w_{2}y_{1}+\frac{w_{2}}{w_{1}}i)\\
=w_{1}^{n}\sum_{k=0}^{n}\binom{n}{k}T_{n-k,q^{w_{1}w_{3}}}(w_{2}-1)w_{3}^{n-k}w_{2}^{k}
\sum_{i=0}^{w_{1}-1}(-1)^{i}q^{w_{2}w_{3}i}E_{k,q^{w_{1}w_{2}}}(w_{3}y_{1}+\frac{w_{3}}{w_{1}}i)\\
=w_{2}^{n}\sum_{k=0}^{n}\binom{n}{k}T_{n-k,q^{w_{1}w_{2}}}(w_{3}-1)w_{1}^{n-k}w_{3}^{k}
\sum_{i=0}^{w_{2}-1}(-1)^{i}q^{w_{1}w_{3}i}E_{k,q^{w_{2}w_{3}}}(w_{1}y_{1}+\frac{w_{1}}{w_{2}}i)\\
=w_{2}^{n}\sum_{k=0}^{n}\binom{n}{k}T_{n-k,q^{w_{2}w_{3}}}(w_{1}-1)w_{3}^{n-k}w_{1}^{k}
\sum_{i=0}^{w_{2}-1}(-1)^{i}q^{w_{1}w_{3}i}E_{k,q^{w_{1}w_{2}}}(w_{3}y_{1}+\frac{w_{3}}{w_{2}}i)\\
=w_{3}^{n}\sum_{k=0}^{n}\binom{n}{k}T_{n-k,q^{w_{1}w_{3}}}(w_{2}-1)w_{1}^{n-k}w_{2}^{k}
\sum_{i=0}^{w_{3}-1}(-1)^{i}q^{w_{1}w_{2}i}E_{k,q^{w_{2}w_{3}}}(w_{1}y_{1}+\frac{w_{1}}{w_{3}}i)\\
=w_{3}^{n}\sum_{k=0}^{n}\binom{n}{k}T_{n-k,q^{w_{2}w_{3}}}(w_{1}-1)w_{2}^{n-k}w_{1}^{k}
\sum_{i=0}^{w_{3}-1}(-1)^{i}q^{w_{1}w_{2}i}E_{k,q^{w_{1}w_{3}}}(w_{2}y_{1}+\frac{w_{2}}{w_{3}}i).\\
\end{split}
\end{equation}
\end{theorem}
Putting $w_{3}=1$ in (\ref{a55}), we obtain the following corollary.
In Section 1, the identities in (\ref{a53}), (\ref{a56}), and
(\ref{a58}) are combined to give those in (\ref{a8})-(\ref{a15}).
\begin{corollary}\label{b12}
  Let $w_{1},w_{2}$ be any odd positive integers.  Then we have:
\begin{equation}\label{a56}
\begin{split}
&w_{1}^{n}\sum_{i=0}^{w_{1}-1}(-1)^{i}q^{w_{2}i}E_{n,q^{w_{1}}}(w_{2}y_{1}+\frac{w_{2}}{w_{1}}i)\\
=&w_{2}^{n}\sum_{i=0}^{w_{2}-1}(-1)^{i}q^{w_{1}i}E_{n,q^{w_{2}}}(w_{1}y_{1}+\frac{w_{1}}{w_{2}}i)\\
=&\sum_{k=0}^{n}\binom{n}{k}E_{k,q^{w_{1}}}(w_{2}y_{1})T_{n-k,q^{w_{2}}}(w_{1}-1)w_{2}^{n-k}w_{1}^{k}\\
=&\sum_{k=0}^{n}\binom{n}{k}E_{k,q^{w_{2}}}(w_{1}y_{1})T_{n-k,q^{w_{1}}}(w_{2}-1)w_{1}^{n-k}w_{2}^{k}\\
=&w_{1}^{n}\sum_{k=0}^{n}\binom{n}{k}T_{n-k,q^{w_{1}}}(w_{2}-1)w_{2}^{k}
\sum_{i=0}^{w_{1}-1}(-1)^{i}q^{w_{2}i}E_{k,q^{w_{1}w_{2}}}(y_{1}+\frac{i}{w_{1}})\\
=&w_{2}^{n}\sum_{k=0}^{n}\binom{n}{k}T_{n-k,q^{w_{2}}}(w_{1}-1)w_{1}^{k}
\sum_{i=0}^{w_{2}-1}(-1)^{i}q^{w_{1}i}E_{k,q^{w_{1}w_{2}}}(y_{1}+\frac{i}{w_{2}}).
\end{split}
\end{equation}
\end{corollary}
Letting further $w_{2}=1$ in (\ref{a56}), we get the following
corollary. This is the identity obtained in
(\ref{a54}) together with the multiplication formula for
 $q$-Euler polynomials
[cf.7,(2.17)].

\begin{corollary}\label{b13}
  Let $w_{1}$ be any odd positive integer.  Then we have:
\begin{equation*}
\begin{split}
E_{n,q}(w_{1}y_{1})&=w_{1}^{n}\sum_{i=0}^{w_{1}-1}(-1)^{i}q^{i}E_{n,q^{w_{1}}}
(y_{1}+\frac{i}{w_{1}})\\
&=\sum_{k=0}^{n}\binom{n}{k}E_{k,q^{w_{1}}}(y_{1})T_{n-k,q}(w_{1}-1)w_{1}^{k}.
\end{split}
\end{equation*}
\end{corollary}
\begin{theorem}\label{b14}
Let $w_{1},w_{2},w_{3}$ be any odd positive integers. Then we have:
\begin{equation}\label{a57}
\begin{split}
&(w_{1}w_{2})^{n}\sum_{i=0}^{w_{1}-1}\sum_{j=0}^{w_{2}-1}(-1)^{i+j}q^{w_{3}(w_{2}i+w_{1}j)}
E_{n,q^{w_{1}w_{2}}}(w_{3}y_{1}+\frac{w_{3}}{w_{1}}i+\frac{w_{3}}{w_{2}}j)\\
=&(w_{2}w_{3})^{n}\sum_{i=0}^{w_{2}-1}\sum_{j=0}^{w_{3}-1}(-1)^{i+j}q^{w_{1}(w_{3}i+w_{2}j)}
E_{n,q^{w_{2}w_{3}}}(w_{1}y_{1}+\frac{w_{1}}{w_{2}}i+\frac{w_{1}}{w_{3}}j)\\
=&(w_{3}w_{1})^{n}\sum_{i=0}^{w_{3}-1}\sum_{j=0}^{w_{1}-1}(-1)^{i+j}q^{w_{2}(w_{1}i+w_{3}j)}
E_{n,q^{w_{1}w_{3}}}(w_{2}y_{1}+\frac{w_{2}}{w_{3}}i+\frac{w_{2}}{w_{1}}j).
\end{split}
\end{equation}
\end{theorem}
Letting $w_{3}=1$ in (\ref{a57}), we have the following corollary.
\begin{corollary}\label{b15}
  Let $w_{1},w_{2}$ be any odd positive integers.  Then we have:
\begin{equation}\label{a58}
\begin{split}
&w_{1}^{n}\sum_{j=0}^{w_{1}-1}(-1)^{j}q^{w_{2}j}E_{n,q^{w_{1}}}(w_{2}y_{1}+\frac{w_{2}}{w_{1}}j)\\
=&w_{2}^{n}\sum_{i=0}^{w_{2}-1}(-1)^{i}q^{w_{1}i}E_{n,q^{w_{2}}}(w_{1}y_{1}+\frac{w_{1}}{w_{2}}i)\\
=&(w_{1}w_{2})^{n}\sum_{i=0}^{w_{1}-1}\sum_{j=0}^{w_{2}-1}(-1)^{i}q^{w_{2}i+w_{1}j}E_{n,q^{w_{1}w_{2}}}
(y_{1}+\frac{i}{w_{1}}+\frac{j}{w_{2}}).
\end{split}
\end{equation}
\end{corollary}
\begin{theorem}\label{b16}
Let  $w_{1},w_{2},w_{3}$ be any positive integers. Then we have:
\begin{equation}\label{a59}
\begin{split}
&\sum_{k+l+m=n}\binom{n}{k,l,m}E_{k,q^{w_{3}}}(w_{1}y)E_{l,q^{w_{1}}}(w_{2}y)
E_{m,q^{w_{2}}}(w_{3}y)w_{3}^{k}w_{1}^{l}w_{2}^{m}\\
=&\sum_{k+l+m=n}\binom{n}{k,l,m}E_{k,q^{w_{2}}}(w_{1}y)E_{l,q^{w_{1}}}(w_{3}y)
E_{m,q^{w_{3}}}(w_{2}y)w_{2}^{k}w_{1}^{l}w_{3}^{m}.
\end{split}
\end{equation}
\end{theorem}
\begin{theorem}\label{b17}
Let $w_{1},w_{2},w_{3}$ be any odd positive integers. Then we have:
\begin{align}
\label{a60}
&\sum_{k+l+m=n}\binom{n}{k,l,m}T_{k,q^{w_{3}}}(w_{1}-1)T_{l,q^{w_{1}}}(w_{2}-1)
T_{m,q^{w_{2}}}(w_{3}-1)w_{3}^{k}w_{1}^{l}w_{2}^{m}\\
\label{a61}
=&\sum_{k+l+m=n}\binom{n}{k,l,m}T_{k,q^{w_{2}}}(w_{1}-1)T_{l,q^{w_{1}}}(w_{3}-1)
T_{m,q^{w_{3}}}(w_{2}-1)w_{2}^{k}w_{1}^{l}w_{3}^{m}.
\end{align}
\\
\end{theorem}
Putting $w_{3}=1$ in (\ref{a60}) and (\ref{a61}), we get the
following corollary.
\begin{corollary}\label{b18}
Let $w_{1},w_{2}$ be any odd positive integers.  Then we have:
\begin{equation*}
\begin{split}
&\sum_{k=0}^{n}\binom{n}{k}T_{k,q^{w_{1}}}(w_{2}-1)T_{n-k,q}(w_{1}-1)w_{1}^{k}\\
=&\sum_{k=0}^{n}\binom{n}{k}T_{k,q^{w_{2}}}(w_{1}-1)T_{n-k,q}(w_{2}-1)w_{2}^{k}.
\end{split}
\end{equation*}
\end{corollary}
\bibliographystyle{amsplain}

\end{document}